\documentclass{jnmp}
\usepackage{amsmath,amssymb,amscd,array}

\let\ssection=\section
\renewcommand{\section}{\setcounter{equation}{0}\ssection}

\def\l{\lambda}

\newcommand{\bbR}{\mathbb{R}}

\newcommand{\Diff}{\mathrm{Diff}}

\newcommand{\End}{\mathrm{End}}
\newcommand{\cF}{{\mathcal{F}}}
\newcommand{\cD}{{\mathcal{D}}}

\newcommand{\Hom}{\mathrm{Hom}}

\newcommand{\PSL}{\mathrm{PSL}}
\newcommand{\SL}{\mathrm{SL}}

\newcommand{\Vect}{\mathrm{Vect}}

\newcommand{\cqfd}{\hspace*{\fill}\rule{3mm}{3mm}}

\def\l{\lambda}

\newtheorem{theorem}{Theorem}
\newtheorem{proposition}{Proposition}
\newtheorem{lemma}{Lemma}

\theoremstyle{definition}

\newtheorem{remark}{Remark}

\begin{document}
\renewcommand{\evenhead}{Sofiane Bouarroudj}
\renewcommand{\oddhead}{On $\mathfrak{sl}(2)$-relative cohomology of the Lie algebra...}

%
\thispagestyle{empty}

\FirstPageHead{14}{1}{2007}{\pageref{firstpage}--\pageref{lastpage}}{Article}

\copyrightnote{2006}{Sofiane Bouarroudj}

\Name{On $\mathfrak{sl}(2)$-relative cohomology of the Lie algebra
of vector fields and differential operators}

\label{firstpage}

\Author{Sofiane Bouarroudj}

\Address{Department of Mathematical Sciences, U.A.E. University,
Faculty of Science, P.O. Box 17551, Al-Ain, U.A.E.

E-mail:bouarroudj.sofiane@uaeu.ac.ae}

\Date{Received June 18, 2006; Accepted in Revised Form 
October 11, 2006}

\begin{abstract}Let $\Vect(\mathbb{R})$ be the Lie algebra of smooth
vector fields on $\mathbb{R}$. The space of symbols
$\mathrm{Pol}(T^* \mathbb{R})$ admits a non-trivial deformation
(given by differential operators on weighted densities) as a
$\Vect(\mathbb{R})$-module that becomes trivial once the action is
restricted to $\mathfrak{sl}(2)\subset \Vect(\mathbb{R}).$ The
deformations of $\mathrm{Pol}(T^* \mathbb{R}),$ which become
trivial once the action is restricted to $\mathfrak{sl}(2)$ and
such that the $\Vect(\mathbb{R})$-action on them is expressed in
terms of differential operators, are classified by the elements of
the weight basis of
$\mathrm{H}^2_{\mathrm{diff}}(\Vect(\mathbb{R}),\mathfrak{sl}(2);
{\cal D}_{\lambda,\mu}),$ where $\mathrm{H}^i_{\mathrm{diff}}$
denotes the differential cohomology (i.e., we consider only
cochains that are given by differential operators) and where
${\cal D}_{\lambda,\mu}=\Hom_{\text{diff}}(\cF_\lambda, \cF_\mu)$
is the space of differential operators acting on weighted
densities. The main result of this paper is computation of this
cohomology. In addition to relative cohomology, we exhibit
2-cocycles spanning $\mathrm{H^2}(\mathfrak{g}; {\cal D}_{\lambda,
\mu})$ for $\mathfrak{g}=\Vect(\mathbb{R})$ and
$\mathfrak{sl}(2).$
\end{abstract}
\section{Introduction}

{\bf Notations.} Let $\Vect(\mathbb{R})$ be the Lie algebra of
smooth vector fields on $\mathbb{R}$. Let ${\cal F}_{\l}$ be the
space of weighted densities of degree $\lambda$ on $\mathbb{R}$,
i.e., the space of sections of the line bundle
$(T^*\mathbb{R})^{\otimes \l}$, so its elements can be represented
as $\phi(x)dx^\lambda$, where $\phi(x)$ is a function and
$dx^\lambda$ is a formal (for a time being) symbol. This space
coincides with the space of vector fields, functions and
differential forms for $\lambda=-1,$ $0$ and $1,$ respectively.
The Lie algebra $\Vect(\mathbb{R})$ acts on $\cal F_\l$ by the Lie
derivative:  we set
\begin{equation}
\label{dens} L_{X}^\l(\phi\,dx^\lambda)=\left (X(\phi) +\l\, \phi
\,\mathrm{div}X\right )dx^\lambda\;\text{for any $X\in
\Vect(\mathbb{R})$ and $\phi\,dx^\lambda \in \cal F_\l$}.
\end{equation}
We denote by ${\cal D}_{\lambda,\nu}$ the space of linear
differential operators that act on the spaces of weighted
densities:
\begin{equation}
\label{Op} A:{\cal F}_{\lambda}\to {\cal F}_{\mu}.
\end{equation}
The Lie algebra $\Vect(\mathbb{R})$ acts on ${\cal
D}_{\lambda,\nu}$ as follows. For any $X\in \Vect(\mathbb{R}),$ we
set (here $L_X^\l$ is the action (\ref{dens})):
\begin{equation}
\label{act} L_X^{\l,\mu}(A)=L_X^\mu \circ A -A \circ L_X^\l.
\end{equation}

\noindent {\bf Motivations.} This work has its genesis in the
study of the $\Vect(\mathbb{R})$-module ${\cal D}_{\lambda, \mu}$.
Duval, Lecomte and Ovsienko showed \cite{do, lo} that this space
cannot be isomorphic, as a $\Vect(\mathbb{R})$-module, to the
corresponding space of symbols of these operators but is its
deformation in the sense of Richardson-Neijenhuis (\cite{nr}). As
is well known, deformation theory of modules is closely related to
the Lie algebra cohomology (\cite{nr}). More precisely, given a
Lie algebra $\mathfrak g$ and a $\mathfrak g$-module $V;$ the {\it
infinitesimal} deformations of the $\mathfrak{g}$-module structure
on $V$, i.e., deformations that are linear in the parameter of
deformation, are described by the elements (up to proportionality)
of $\mathrm{H}^1(\mathfrak g; \End(V))$. The obstructions to
extension of any infinitesimal deformation to a formal one are
similarly described by $\mathrm{H}^2(\mathfrak g; \End(V))$.
Computation of $\mathrm{H}^1$ in our situation (with
$\mathfrak{g}=\Vect(\mathbb{R})$ and ${\cal D}_{\lambda, \mu}$
instead of $\End(V)$) was carried out by Feigin and Fuchs
\cite{ff}. Ovsienko and I computed the corresponding
$\mathfrak{sl}(2)$-relative cohomology (see \cite{bo}). Gordan's
classification of bilinear differential operators on weighted
densities \cite{g} played a central role in our computation.
Later, a generalization to multi-dimensional manifolds has been
carried out by Lecomte and Ovsienko in \cite{lo}; for further
results, see \cite{b2}. Note that the $\mathfrak{sl}(2)$-relative
cohomology measures infinitesimal deformations that become trivial
once the action is restricted to $\mathfrak{sl}(2).$ This is
actually the case for the space of differential operators since,
as $\mathfrak{sl}(2)$-module, it is isomorphic to the space of
symbols for generic $\lambda$ and $\mu$ (cf. \cite{gar}). Let
$\mathrm{H}^i_{\mathrm{diff}}$ be the differential cohomology
(i.e., we consider only cochains that are given by differential
operators). Recently I realized that a description of (here
$\Vect_{\mathrm P}(\mathbb{R})$ is the Lie algebra of polynomial
vector fields)
\begin{equation}
\label{intr}
\mathrm{H}^2_{\mathrm{diff}}(\Vect_{\mathrm{P}}(\mathbb{R}); {\cal
D}_{\lambda, \mu})
\end{equation}
can be deduced from the work by Feigin and Fuchs \cite{ff}.
Feigin-Fuchs gave details of computation of $
\mathrm{H}^1_{\mathrm{diff}}(\Vect_{\mathrm{P}}(\mathbb{R}); {\cal
D}_{\lambda, \mu})$ but not of higher cohomology and no explicit
2-cocycles were provided. The $\mathfrak{sl}(2)$-relative
cohomology cannot, however, be deduced from their computation.
Several authors (see, e.g., \cite{los, tsu}) have also studied
${\mathrm H}^i(\Vect(\mathbb{R}); {\cal A})$ for an arbitrary
$\Vect(\mathbb{R})$-module ${\cal A}.$ But it is not easy to get a
description of the cohomology (\ref{intr}) nor the
$\mathfrak{sl}(2)$-relative cohomology from their results. Our
main result is computation of the $\mathfrak{sl}(2)$-relative
cohomology and explicit expressions of 2-cocycles that span
(\ref{intr}). This work is the first step towards the study of
formal deformations of symbols.

For investigation of all deformations of symbols in case of
$\mathbb{R}^n$ for $n>1$, see \cite{aalo}. The authors use the
Neijenhuis-Richardson product to prove the existence of cocycles
but do not compute any cohomology. The cohomology similar to
(\ref{intr}) with $\mathbb{R}^n$ instead of $\mathbb{R}$ is still
out of reach for $n>1$.
\section{Basic definitions}
Consider the standard (local) action of $\SL(2)$ on $\mathbb{R}$
by linear-fractional transformations. Although the action is
local, it generates global vector fields
$$
\frac{d}{dx},\quad x\frac{d}{dx},\quad x^2\frac{d}{dx},
$$ that form a Lie subalgebra of $\Vect(\mathbb{R})$ isomorphic
to the Lie algebra $\mathfrak{sl}(2)$ (cf. \cite{olv}). This
realization of $\mathfrak{sl}(2)$ is understood throughout this
paper.

\subsection{The Gelfand-Fuchs cocycle}
We need to introduce the following cocycle (of Gelfand-Fuchs):
\begin{equation}
\label{gelfuk} \omega(X,Y)=\left |
\begin{array}{cc}
f' & g''\\
f' & g''
\end{array}
\right |\, dx\quad \mbox{ for } X=f\frac{d}{dx}, Y=g\frac{d}{dx}.
\end{equation}
Here $\omega$ is a cohomology class in ${\mathrm
H}^2(\Vect(\mathbb{R}), \cF_1).$ Related is the element of
$\mathrm{H^2}(\mathrm{Vect}(S^1)),$ the 2-cocycle on
$\mathrm{Vect}(S^1)$ given by the formula (see \cite{gf}):
$$
\int_{S^1}\,\omega(X,Y).
$$
This 2-cocycle generates the central extension of $\Vect(S^1)$
called the {\it Virasoro} algebra.
\section{The $\mathfrak{sl}(2)$-relative cohomology of
$\Vect(\mathbb{R})$ acting on ${\cal D}_{\lambda,\mu}$}
The following steps to compute the relative cohomology has
intensively been used in \cite{b1,b2,bo,lo}. First, we classify
$\mathfrak{sl}(2)$-invariant differential operators, then we
isolate among them those that are 2-cocycles. To do that, we need
the following Lemma.
\begin{lemma}
\label{har} Any 2-cocycle vanishing on the Lie subalgebra
$\mathfrak{sl}(2)$ of $\Vect(\mathbb{R})$ is
$\mathfrak{sl}(2)$-invariant.
\end{lemma}
\begin{proof} The 2-cocycle condition reads as follows:
$$
\begin{array}{l}
c([X,Y],Z,\phi\,dx^\lambda)-L_X^{\lambda,\mu}\,c(Y,Z,\phi\,dx^\lambda)
+ \circlearrowleft (X,Y,Z)=0
\end{array}
$$
for every  $X, Y, Z\in \Vect(\mathbb{R})$ and $\phi\,dx^\lambda\in
\cF_\lambda$, where $ \circlearrowleft (X,Y,Z)$ denotes the
summands obtained from the two written ones by the cyclic
permutation of the symbols $X, Y, Z$. Now, if $X\in
\mathfrak{sl}(2),$ then the equation above becomes
$$
c([X,Y],Z,\phi\,dx^\lambda)-
c([X,Z],Y,\phi\,dx^\lambda)=L_X^{\lambda,\mu}\,c(Y,Z,\phi\,dx^\lambda).
$$
This condition is nothing but the invariance property.
\end{proof}
\subsection{$\mathfrak{sl}(2)$-invariant differential operators}
\label{emi}
As our 2-cocycles vanish on $\mathfrak{sl}(2),$ we will
investigate $\mathfrak{sl}(2)$-invariant bilinear differential
operators that vanish on $\mathfrak{sl}(2).$
\begin{proposition}
\label{gh} The space of skew-symmetric bilinear differential
operators $ \Vect(\mathbb{R})\wedge \Vect(\mathbb{R})\rightarrow
\cD_{\lambda, \mu}, $ which are $\mathfrak{sl}(2)$-invariant and
vanish on $\mathfrak{sl}(2),$ is as follows:
\begin{enumerate}
\item It is $\frac{1}{2}(k-3)$-dimensional if $\mu-\lambda=k$ and
$k$ is odd. \item It is $\frac{1}{2}(k-4)$-dimensional if
$\mu-\lambda=k$ and $k$ is even. \item It is $0$-dimensional,
otherwise.
\end{enumerate}
\end{proposition}
\begin{proof} The generic form of any such a differential operator is
(here $X=f \frac{d}{dx},Y=g \frac{d}{dx}\in \Vect(\mathbb{R})$ and
$\phi\,dx^\lambda\in \cF_\lambda$):
$$
c(X,Y,\phi\,dx^\lambda)=\sum_{i+j+l\leq k}\,
c_{i,j}\,f^{(i)}\,g^{(j)}\,\phi^{(l)}dx^{\mu},
$$
where $c_{i,j}=-c_{j,i}$ and $f^{(i)}$ stands for $\frac{d^{i}
f}{dx^{i}}.$

The invariance property with respect to the vector field
$X=x\frac{d}{dx}$ with arbitrary $Y$ and $Z$ implies that
$c_{i,j}'=0$ and $\mu=\lambda+i+j+l.$ Therefore $c_{i,j}$ are
constants. Now, the invariance property with respect $X=x^2
\frac{d}{dx}$ with arbitrary $Y$ and $Z$ is equivalent to the
system (where $2<\beta<\gamma<k$):
\begin{equation}
\label{par} \small (\beta+1)(\beta-2)\, c_{\beta+1,\gamma}-
(\gamma+1)(\gamma-2)\, c_{\gamma+1,\beta}+
(k+2-\beta-\gamma)(k+1-\beta-\gamma+2\lambda)\,
c_{\beta,\gamma}=0.
\end{equation}
For $\beta=3,$ the equation (\ref{par}) implies that all the
constants $c_{t,3}$ can be determined {\it uniquely} in terms of
$c_{4,3}$ and $c_{4,s}.$ More precisely,
$$
c_{\gamma+1,3}=\frac{4\, c_{4,\gamma}+(k-1-\gamma)\left (
k-2-\gamma+2\, \lambda\right )\,
c_{3,\gamma}}{(\gamma+1)(\gamma-2)}.
$$
For $\beta=4$ and $\gamma=5,$ and from the system (\ref{par}), we
have
$$
c_{6,4}=\frac{1}{12}\, (k-7)\left (k-8+2\,\lambda\, \right )\,
c_{4,5}.
$$
Thus the constant $c_{6,4}$ is determined. But for $\beta=4$ and
$\gamma>5,$ the system (\ref{par}) implies that
$$
c_{5,\gamma}=\frac{1}{10}(\gamma+1)(\gamma-2)\,
c_{\gamma+1,4}-\frac{1}{10}(k-\gamma-2)(k-\gamma-3 +2\, \lambda)\,
c_{4,\gamma}.
$$
Therefore all $c_{5,\gamma}$ can be determined for any $\gamma\geq
6.$

By continuing this procedure we see that $c_{6,\gamma},
c_{7,\gamma},\ldots$ can be determined as well as $c_{4,\gamma}$
for $\gamma$ even. Finally, we have proved that the space of
$\mathfrak{sl}(2)$-invariant operators is as follows:

(i) for $k$ even, it is generated by $c_{4,3}, c_{4,5}, c_{4,7},
\ldots, c_{4,k-3}.$ The space of solution is
$\frac{1}{2}(k-4)$-dimensional.

(ii) for $k$ odd, it is generated by $c_{4,3}, c_{4,5}, c_{4,7},
\ldots, c_{4,k-2}.$ The space of solution is
$\frac{1}{2}(k-3)$-dimensional.
\end{proof}
\subsection{The $\mathfrak{sl}(2)$-relative cohomology of
$\Vect(\mathbb{R})$}
\begin{theorem}
\label{main} We have
\begin{equation} \nonumber {\mathrm
H}^2_{\mathrm{diff}}(\Vect(\mathbb{\mathbb{R}}),\mathfrak{sl}(2);
\cD_{\lambda,\mu})= \left\{
\begin{array}{ll}
\bbR& \mbox{if } \left \{
\begin{array}{l}
(\lambda,\mu)=(0,5), (-2,3), (-4,1), (-\frac{5}{2}, \frac{7}{2})
\mbox{ or } \\[2mm]
\quad \quad \quad \,\,\,\, \left (-\frac{5}{2}\pm
\frac{\sqrt{19}}{2},
\frac{7}{2}\pm \frac{\sqrt{19}}{2}\right ),\\[2mm]
\mu-\lambda=7,8,9,10,11 \mbox{ for all }
\lambda\not =\frac{1-k}{2}, \\[2mm]
\mu-\lambda=k=12,13,14 \mbox{ and } \lambda=
\frac{1-k}{2}\pm \frac{\sqrt{12k-23}}{2},\\[2mm]
\mu-\lambda=k=15 \mbox{ and } \lambda=\frac{1-k}{2},\\[2mm]
\end{array}
\right.
\\
\mathbb{R}^2& \mbox{if }
\begin{array}{l}\,\,\mu-\lambda=k=7,\ldots, 14\mbox{ for }
\lambda=\frac{1-k}{2},
\\
\end{array}
\\[2mm]
0& \mbox{otherwise. }
\\
\end{array}
\right.
\end{equation}
\end{theorem}
\begin{remark}{\rm $ {\mathrm H}^1
_{\mathrm{diff}}(\Vect(\mathbb{\mathbb{R}}),\mathfrak{sl}(2);
\cD_{\lambda,\mu})$ has been computed in \cite{bo}.}
\end{remark}
\section{Proof of Theorem \ref{main}}
Every 2-cocycle on $\Vect(\mathbb{R})$ retains the following
general form (here $X=f \frac{d}{dx},Y=g \frac{d}{dx}\in
\Vect(\mathbb{R})$ and $\phi\,dx^\lambda\in \cF_\lambda$):
\begin{equation}
\label{sam} c(X,Y,\phi\,dx^\lambda)=\sum_{i+j+l\leq k}\,
c_{i,j}\,f^{(i)}\,g^{(j)}\,\phi^{(l)}dx^{\mu},
\end{equation}
where $c_{i,j}=-c_{j,i}.$ Since this 2-cocycle vanishes on
$\mathfrak{sl}(2)$, Lemma \ref{har} implies that this 2-cocycle is
$\mathfrak{sl}(2)$-invariant. Therefore all $c_{i,j}$ are zero and
$i+j+l=\mu-\lambda.$ The last statement means that the 2-cocycle
(\ref{sam}) is homogenous. Besides, we have
$c_{0,j}=c_{1,j}=c_{2,j}=0.$

Before starting with the proof proper, we explain our strategy. This
method has already been used in \cite{b1}. First, we investigate
operators that belong to
$Z^2(\Vect(\mathbb{\mathbb{R}}),\mathfrak{sl}(2);
\cD_{\lambda,\mu}).$ The 2-cocycle condition imposes conditions on
the constants $c_{i,j}$: we get a linear system for $c_{i,j}.$
Second, taking into account these conditions, we eliminate all
constants underlying coboundaries. Gluing these bits of information
together we deduce that $\dim {\mathrm H}^2$ is equal to the number
of independent constants $c_{i,j}$ remaining in the expression of
the 2-cocycle (\ref{sam}).

\begin{proposition}{\em (\cite{g})}
\label{gor} There exist $\mathfrak{sl}(2)$-invariant bilinear
differential operators $J_k^{\tau,\lambda}: \cF_\tau\otimes
\cF_\lambda \rightarrow \cF_{\tau+\lambda+k}$ given by:
\begin{equation}
\label{hmida} J_k^{\tau,\lambda}(\varphi\,dx^\tau,
\phi\,dx^\lambda) =\sum_{i+j=k}
\gamma_{i,j}\,\varphi^{(i)}\,\phi^{(j)}\, dx^{\tau+\lambda+k},
\end{equation}
where the constants $\gamma_{i,j}$ satisfy
\begin{equation}
\label{ak}
(i+1)(i+2\tau)\,\gamma_{i+1,j}+(j+1)(j+2\lambda)\,\gamma_{i,j+1}=0.
\end{equation}
\end{proposition}
\begin{remark}{\rm The operators (\ref{hmida}) are called {\it
transvectants}. Amazingly, they appear in many contexts, especially
in the computation of cohomology (cf. \cite{b1, bo}). We refer to
\cite{ot} for their history.}
\end{remark}
Now we will study properties of the coboundaries. Let
$B:\Vect(\mathbb{R})\rightarrow \cD_{\lambda,\mu}$ be an operator
defined by (for any $X=f\frac{d}{dx}\in \Vect(\mathbb{R})$ and
$\phi\,dx^\lambda \in \cF_\lambda$):
$$
B(X,\phi\,dx^\lambda)=\sum_{i+j=k+1}\,\gamma_{i,j}\, f^{(i)}\,
\phi^{(j)}\, dx^{\lambda+k}.
$$
\begin{proposition}\label{maa}
Every coboundary $\delta  (B)\in B^2(\Vect(\mathbb{R}),
\mathfrak{sl}(2); \cD_{\lambda,\mu})$ possesses the following
properties. The operator $B$ coincides (up to a nonzero factor)
with the transvectant $J_{k+1}^{-1,\lambda}$, where
$\gamma_{0,k+1}=\gamma_{1,k}=\gamma_{2,k-1}=0.$ In addition (here
$X=f\frac{d}{dx}\in \Vect(\mathbb{R})$ and $\phi\,dx^\lambda \in
\cF_\lambda$),
\begin{equation}
\label{nedj}\delta (B)(X,Y,\phi\,dx^\lambda)=
\sum_{i+j+l=k+2}\,\beta_{i,j}\,f^{(i)}\,g^{(j)}\,\phi^{(l)}\,
dx^{\lambda+k},
\end{equation}
where
$$
\beta_{0,j}=\beta_{1,j}=\beta_{2,j}=0,
$$
and
\small{
$$
\begin{array}{lcl}
\beta_{3,4}&=&-\frac{1}{24}\binom{k-2}{3}\left(k^2 + 4 (\lambda
-1)\lambda + k (4\lambda-5)\right )( k -1+ 2 \lambda
)\gamma_{3,k-2}
\\[3mm]
\beta_{4,5}&=&-\frac{1}{480}\binom{k-2}{5} (k-3+ 2\lambda )(k^3 +
4(\lambda-1) \lambda (2\lambda-19) + 3 k^2(2
\lambda-7 )+2k(49+6(\lambda-7 )\lambda))\\[3mm]
&&\times (k -1+ 2 \lambda )\gamma_{3,k-3}.
\end{array}
$$
}
\end{proposition}
\begin{proof}
From the very definition of coboundaries, we have (for any $X,
Y\in \Vect(\mathbb{R})$ and $\phi\,dx^\lambda\in \cF_\lambda$):
$$
\delta(B)(X,Y,\phi\,dx^\lambda)=B([X,Y],\phi\,dx^\lambda)-L_X\,
B(Y,\phi\,dx^\lambda)+L_Y\, B(X,\phi\,dx^\lambda).
$$
The coboundary above vanishes on the Lie algebra
$\mathfrak{sl}(2).$ It means that if $X\in \mathfrak{sl}(2)$, we
have
$$
B([X,Y],\phi\,dx^\lambda)=L_X\, B(Y,\phi\,dx^\lambda)-L_Y\,
B(X,\phi\,dx^\lambda).
$$
Hence, the operator $B$ is $\mathfrak{sl}(2)$-invariant; therefore
it coincides with the transvectants. The conditions
$\gamma_{0,k+1}=\gamma_{1,k}=\gamma_{2,k-1}=0$ come from the fact
that the operator $B$ vanishes on $\mathfrak{sl}(2).$ Now, the
conditions $\beta_{0,j}=\beta_{1,j}=\beta_{2,j}=0$ are
consequences of $\mathfrak{sl}(2)$-invariance, while the values of
$\beta_{3,4}$ and $\beta_{4,5}$ follow by a direct computation.
\end{proof}
\subsection{The case where $\mu-\lambda=5$}
In this case, the 2-cocycle has the form
\begin{equation}
\label{or} c(X,Y,\phi\,dx^\lambda)= \left |
\begin{array}{ll}
f^{(3)}& g^{(3)}\\
f^{(4)}& g^{(4)} \end{array} \right |\, \phi\,dx^{\lambda+5}\quad
\mbox{for } X=f\frac{d}{dx}, Y=g\frac{d}{dx}.
\end{equation}
The 2-cocycle condition is always satisfied. On the other hand,
the coboundary (\ref{nedj}) takes the form
$$
\frac{1}{3}\lambda \,\left( 2 + \lambda  \right) \,\left( 4 +
\lambda \right) \,
     {{\gamma }_{3, k-2}}\,\left( g^{(3)}\,f^{(4)} -
       f^{(3)}\,g^{(4)} \right)\,\phi\,dx^{\lambda+5}.
$$
This coboundary coincides with the 2-cocycle (\ref{or}) except for
$\lambda =0, -2$ or $-4.$ Therefore the cohomology in Theorem
\ref{main} is trivial except for $\lambda=0,-2$ or $-4.$
\subsection{The case where $\mu-\lambda=6$}
The 2-cocycle has the form
$$
c(X,Y,\phi\,dx^\lambda)=\left ( \left |
\begin{array}{ll}
f^{(3)}& g^{(3)}\\
f^{(4)}& g^{(4)} \end{array} \right |\right .\,
\phi'-\frac{\lambda}{5}\, \left .\left |
\begin{array}{ll}
f^{(3)}& g^{(3)}\\
f^{(5)}& g^{(5)} \end{array} \right |\phi  \right )
dx^{\lambda+6}\quad \mbox{ for } X=f\frac{d}{dx}, Y=g\frac{d}{dx}.
$$
On the other hand, the coboundary (\ref{nedj}) takes the form
$$
\begin{array}{l}
\frac{1}{3}\left( 5 + 2\,\lambda  \right) \,  \left( 3 +
2\,\lambda \,\left( 5 + \lambda  \right)  \right) \,{{\gamma
}_{3,k-2}}\, \,\left( g^{(3)}\,f^{(4)} - f^{(3)}\,g^{(4)} \right)
\phi'\,dx^{\lambda+6}\\[2mm]
- \frac{1}{15}\lambda \,\left( 5 + 2\,\lambda  \right) \,
      \left( 3 + 2\,\lambda \,\left( 5 + \lambda  \right)
       \right) \,{{\gamma }_{3,k-2}}\,
\,\left( g^{(3)}\,f^{(5)} - f^{(3)}\,g^{(5)}
\right)\phi\,dx^{\lambda+6}.
\end{array}
$$
This coboundary coincides with our 2-cocycle except when
$\lambda=-\frac{5}{2}$ or $\lambda$ is a solution to
$3+2\lambda(5+\lambda)=0.$
\subsection{The case where $\mu-\lambda\geq 7$}
\label{hop}
In this case, the 2-cocycle condition is equivalent to the system
(where $2\leq\alpha <\beta<\gamma$):
\begin{equation}
\label{ref}
\begin{array}{l} \left
(\binom{\alpha+\beta-1}{\alpha}-\binom{\alpha+\beta-1}{\alpha-1}\right
) c_{\alpha+\beta-1, \gamma}- \left
(\binom{\alpha+\gamma-1}{\alpha}-\binom{\alpha+\gamma-1}{\alpha-1}\right
) c_{\alpha+\gamma-1, \beta}\\[2mm]
+\left
(\binom{\beta+\gamma-1}{\beta}-\binom{\beta+\gamma-1}{\beta-1}\right
) c_{\beta+\gamma-1, \alpha} +\left
(\binom{k+2-\beta-\gamma}{\alpha}+\lambda
\binom{k+2-\beta-\gamma}{\alpha-1}\right
) c_{\beta, \gamma}\\[2mm]
-\left (\binom{k+2-\alpha-\gamma}{\beta}+\lambda
\binom{k+2-\alpha-\gamma}{\beta-1}\right ) c_{\alpha, \gamma} +
\left (\binom{k+2-\alpha-\beta}{\gamma}+\lambda
\binom{k+2-\alpha-\beta}{\gamma-1}\right ) c_{\alpha, \beta}=0.
\end{array}
\end{equation}
This system can be deduced by a simple computation. Of course,
such a system has at least one solution in which the solutions
$c_{i,j}$ are just the coefficients $\beta_{i,j}$ of the
coboundaries (\ref{nedj}).

\subsubsection{The case where $\mu-\lambda=7,8,9,10,11$}
\label{coc}
Let us show that the solutions to the system (\ref{ref}) are
expressed in terms of $c_{3,4}$ and $c_{4,5}$.

In the case $\alpha=2,$ the system (\ref{ref}) has been studied in
Section \ref{emi}; its study corresponds to the investigation of
$\mathfrak{sl}(2)$-invariant differential operators. We have seen
that all the constants $c_{i,j}$ can be expressed in terms of
$c_{3,4}, c_{5,4}, c_{7,4}, c_{9,4},\ldots.$

{\bf For $k=7.$} According to Proposition \ref{gh}, the space of
solutions is generated by $c_{3,4}$ and $c_{4,5}$. Note that the
coefficients $c_{4,i},$ where $i\geq 6$, are zero. The following
coefficients can be deduced from the system (\ref{ref}):
\begin{equation}
\small \label{que1}
\begin{array}{lcllcl}
c_{3,5}&=&\frac{1}{10} (5-k)(k-6+2\, \lambda)\, c_{3,4},&
c_{3,6}&=& \frac{1}{18} ((6-k)(k-7+2\, \lambda)\, c_{3,5}-4\,
c_{4,5}).
\end{array}
\end{equation}

{\bf For $k=8$.} According to Proposition \ref{gh}, the space of
solutions is generated by $c_{3,4}$ and $c_{4,5}$. Moreover, the
coefficients  $c_{4,i},$ where $i\geq 7,$ are zero. The solutions
to (\ref{ref}) are given by (\ref{que1}) together with
\begin{equation}
\small
\label{que2}
\begin{array}{lcllcl}
c_{3,7}&=&\frac{1}{28}((7-k)(k+2(\lambda-4))\, c_{3,6}- 4\,
c_{4,6}),& c_{4,6}&=&\frac{1}{18}(k-7)(k-8+2\, \lambda)\, c_{4,5}.
\end{array}
\end{equation}
Now for $k=9,10$ and $11$ we have to deal with the system
(\ref{ref}) for $\alpha=3:$
\begin{equation}
\nonumber
\begin{array}{l}
\left (\binom{\beta+2}{3}-\binom{\beta+2}{2}\right )\,
c_{\beta+2,\gamma}- \left
(\binom{\gamma+2}{3}-\binom{\gamma+2}{2}\right )\,
c_{\gamma+2,\beta} + \left
(\binom{\gamma+\beta-1}{\beta}-\binom{\gamma+\beta-1}{\beta-1}\right
)\,
c_{\gamma+\beta-1,3}\\[2mm]+ \left (\binom{k+2-\beta-\gamma}{3}+\lambda
\binom{k+2-\beta-\gamma}{2}\right )\, c_{\beta,\gamma} -\left
(\binom{k-1-\gamma}{\beta}+\lambda
\binom{k-1-\gamma}{\beta-1}\right )\, c_{3,\gamma}\\[2mm]+ \left
(\binom{k-1-\beta}{\gamma}+\lambda \binom{k-1-\beta}{2}\right )\,
c_{3,\beta} =0.
\end{array}
\end{equation}
For $\beta=4$ and $\gamma=5,$ the coefficient $c_{4,7}$ is given
by
\begin{equation} \small{
\begin{array}{rl}
c_{4,7}=\frac{1}{105840}\binom{k-7}{2}&\left ( \binom{k-5}{2}(2\,
\lambda+k-3)(-288+k(194+k(k-27))+268\,\lambda+6(k-18)k\lambda \right.\\
& \left. +12(k-9)\lambda^2+8\lambda^3)\,c_{3,4}-80
c_{4,5}(279+2k^2+\lambda(8\lambda-113)+k(8\lambda-49)))  \right ).
\end{array}}
\label{que1p}
\end{equation}
We continue like this until we determine all the constants
$c_{4,k-3}$ for $k$ even and $c_{4,k-2}$ for $k$ is odd. Therefore
the system (\ref{ref}) admits solutions generated by $c_{3,4}$ and
$c_{4,5}.$ Let us give explicitly these solutions. \\

{\bf For $k=9$.} The coefficients are given by (\ref{que1}),
(\ref{que2}), (\ref{que1p}) together with
\begin{equation}
\small{ \label{que3}
\begin{array}{lcllcl}
c_{3,8}&\!=\!& \frac{1}{40} ((8-k)(k-9+2\,\lambda)\, c_{3,7}-4\,
c_{4,7}),& c_{5,6}&\!=\!&
\frac{1}{45}\binom{k-8}{2}\,\binom{k+2\lambda-7}{2}\,
c_{3,4}-\frac{14}{5}c_{4,7}.
\end{array}
}
\end{equation}

{\bf For $k=10$.} The coefficients are given by (\ref{que1}),
(\ref{que2}), (\ref{que1p}), (\ref{que3}) together with
\begin{equation}
\small
\label{que4}
\begin{array}{lcllcl}
\small
c_{3,9}&\!\!\!\!\!=\!\!\!\!\!&\frac{1}{54}((9-k)(k+2\lambda-10)\,
c_{3,8} -4\,
c_{4,8}),&c_{5,7}&\!\!\!\!\!=\!\!\!\!\!&
\frac{1}{10}(9-k)(k-10+2\lambda)\,c_{4,7}-4\,
c_{4,8},
\end{array}
\end{equation}
and
\begin{equation}
\label{que5} \small
\begin{array}{ccl}
c_{4,8}&=&\frac{1}{20160}(9-k)(k+2(\lambda-5))\times \\[2mm]
&&\left ((k-8)(k-7)(k+2(\lambda-4))(k-9+2\lambda)\,c_{3,4}+10008\,
c_{4,7}\right ).
\end{array}
\end{equation}

{\bf For $k=11.$} The coefficients are given by (\ref{que1}),
(\ref{que2}), (\ref{que1p}), (\ref{que3}), (\ref{que4}),
(\ref{que5}) together with
\begin{equation}\nonumber
\small
\begin{array}{lcllcl}
c_{3,10}&\!\!=\!\!& \frac{1}{70}((10-k)(k-11+2\,\lambda)\, c_{3,9}
-4\, c_{4,9}),&
c_{5,8}&\!\!=\!\!&\frac{1}{10}(10-k)(k-11+2\lambda)\,
c_{4,8}-\frac{27}{5}\, c_{4,9},
\end{array}
\end{equation}
and
\begin{equation}
\small
\begin{array}{lcl}
 c_{6,7}&=&\frac{1}{45360}\left( \left(k-10 \right) \,\left(k -9
\right) \,
     \left( k + 2\,\left( \lambda-5  \right)  \right) \,
     \left(k + 2\,\lambda -11 \right) \times \right.\\[2mm]
     &&\left .
    ( (k -8) \,(k-7)\, ( k + 2\, \lambda -8 )\,
    ( k + 2\,\lambda -9 ) \,c_{3,4}+ 756\,c_{7,4})\right ) +
        12\,c_{4,9}.\\

        \end{array}
\end{equation}
The explicit value of $c_{4,9}$ is too long; hereafter we omit
such expressions obtained with the help of {\it Mathematica}.

We have just proved that the coefficients of every 2-cocycle is
expressed in terms of the two constants $c_{3,4}$ and $c_{4,5}.$
But this general formula may contain coboundaries. We explain how
the coboundaries can be removed. Consider any coboundary given as
in (\ref{nedj}). We discuss the following cases:

1) $\lambda=\frac{1-k}{2}.$ Then the constant $\beta_{3,4}$ and
$\beta_{4,5}$ vanish simultaneously. Hence the constants $c_{4,5}$
and $c_{3,4}$ cannot be eliminated by adding the coboundary
(\ref{nedj}). It follows that the coefficients of the 2-cocycle
are generated by $c_{3,4}$ and $c_{4,5}.$ Therefore the cohomology
is two-dimensional. The 2-cocycles are given explicitly by the
constants (\ref{que1}), (\ref{que2}), (\ref{que1p}), (\ref{que3}),
(\ref{que4}), (\ref{que5}) by taking $c_{3,4}=1$ and $c_{4,5}=0$
then by taking $c_{3,4}=0$ and $c_{4,5}=1.$

2) $\lambda=\frac{1}{2}(1 - k \pm \sqrt{1 + 3 k}).$ Then the
constant $c_{4,5}$ can be eliminated by adding the coboundary
(\ref{nedj}). On the other hand, the constant $c_{3,4}$ cannot be
eliminated because $\beta_{3,4}=0.$ It follows that the
coefficients of the 2-cocycle are generated by $c_{3,4}.$
Therefore the cohomology is one-dimensional. The 2-cocycle is
given explicitly by the constants (\ref{que1}), (\ref{que2}),
(\ref{que1p}), (\ref{que3}), (\ref{que4}), (\ref{que5}) upon
taking $c_{3,4}=1$ and $c_{4,5}=0.$

3) $\lambda=\frac{3-k}{2}.$ First, we observe that there is no
common solutions for $\lambda$ in 2) and 3) except for $\lambda=1$
and $k=1;$ or $\lambda=-1$ and $ k=1.$ But these cases are not
taken into consideration because $k\geq 7.$ The constant $c_{3,4}$
can be eliminated by adding the coboundary (\ref{nedj}). On the
other hand, the constant $c_{4,5}$ cannot be eliminated because
$\beta_{4,5}=0.$ It follows that the coefficients of the 2-cocycle
are generated by $c_{4,5}.$ Therefore the cohomology is
one-dimensional. The 2-cocycle is given by the constants
(\ref{que1}), (\ref{que2}), (\ref{que1p}), (\ref{que3}),
(\ref{que4}), (\ref{que5}) upon taking $c_{3,4}=0$ and
$c_{4,5}=1.$

4) $\lambda$ is a solution to the equation
$$
k^3 + 4(\lambda-1) \lambda (2\lambda-19) + 3 k^2(2 \lambda-7
)+2k(49+6(\lambda-7 )\lambda)=0.
$$
In this case, $c_{3,4}$ can be eliminated by adding the coboundary
(\ref{nedj}). On the other hand, the constant $c_{4,5}$ cannot be
eliminated as $\beta_{4,5}=0.$ It follows that the coefficients of
the 2-cocycle are generated by $c_{4,5}.$ Therefore the cohomology
is one-dimensional. The coefficients of the 2-cocycle are given by
constants as above upon taking $c_{3,4}=0$ and $c_{4,5}=1.$

5) $\lambda$ is not like 1)--4). In this case, whatever the weight
$\lambda$ is, one of the constant $c_{3,4}$ or $c_{4,5}$ can be
eliminated by adding the coboundary. It follows that the
cohomology is one-dimensional. The coefficients of the 2-cocycle
are given by the constants as above upon taking, for instance,
$c_{3,4}=1$ and $c_{4,5}=0.$
\subsubsection{The case where $\mu-\lambda=12, 13, 14$}
Let us prove that the system (\ref{ref}) has solutions that can be
expressed in terms of one parameter if $\lambda$ is generic, and
in terms of two parameters for particular values of $\lambda.$ But
we have already seen in the previous section that all the
solutions can be expressed in terms of $c_{3,4}$ and $c_{3,5}.$ As
$k\geq 12,$ we are required to study (\ref{ref}) for $\alpha=4.$
For $\alpha=4, \beta=5$ and $\gamma=6,$ the system has one more
equation
\begin{equation}\label{bis}
\small
\begin{array}{l}
\binom{k-7}{5}\,(2\lambda+k-1) \,\left( k^2 + 2\,k\,\left(
2\,\lambda -7  \right)  +
    4\,\left( 6 + \left( -1 + \lambda  \right) \,\lambda  \right)  \right)\times
    \\[1mm]
\left [ \left(k-5 \right) \,k\,
     \left( \left(k-14 \right) \,\left( k-7 \right) \,
        \left( k -6\right) \,\left( k-3 \right) \,c_{3,4} - 400\,c_{4,5} \right)
        \right.\\[1mm]
        \left.
    +4\,\left( \left( k-6 \right) \,\left(  k-5 \right) \,
        \left( -57 + k\,\left( 131 + 2\,\left( -18 + k \right) \,k \right)
              \right) \,c_{3,4} - 400\,\left(  k-1 \right) \,c_{4,5} \right)
              \,\lambda\right .
        \\[1mm]
        \left .
     + 4\,\left( \left(  k-6 \right) \,\left(  k-5 \right) \,
        \left( 101 + 6\,\left( k-12 \right) \,k \right) \,c_{3,4}- 400\,c_{4,5}
       \right) \,{\lambda }^2 \right.\\[1mm]
       \left .
    +32\,{\left(  k-6 \right) }^2\,\left( k-5 \right) \,c_{3,4}\,{\lambda }^3 +
    16\,\left( k-6 \right) \,\left(  k -5\right) \,c_{3,4}\,{\lambda }^4
    \right ]=0.
     \end{array}
\end{equation}
We have three cases:

1) If $\lambda\neq\frac{1}{2}(1-k)$ or $ \frac{1}{2}(1-k\pm
\sqrt{12k-23}),$ then from Eq. (\ref{bis}) the constant $c_{4,5}$
can be expressed in terms of $c_{3,4}.$ Here we have two subcases:

1.1) If $\lambda=\frac{1}{2}(1-k\pm \sqrt{1+3k}),$ then Eq.
(\ref{bis}) implies that $c_{3,4}=0.$ The constant $c_{4,5}$ can
be eliminated by adding the coboundary (\ref{nedj}) for a suitable
$\gamma_{3,k-2}.$ Therefore the cohomology is zero.

1.2) If $\lambda\neq\frac{1}{2}(1-k\pm \sqrt{1+3k}),$ then Eq.
(\ref{bis}) implies that $c_{4,5}$ can be determined in terms of
$c_{3,4}.$ We omit here the explicit expression because it is too
long.

The constant $c_{3,4}$ can be eliminated upon adding the
coboundary (\ref{nedj}) for a suitable $\gamma_{3,k-2}.$ Therefore
the cohomology is zero.

\noindent 2) If $\lambda=\frac{1}{2}(1-k\pm \sqrt{12k-23})$, then
the system (\ref{ref}) has solutions that still depend on
$c_{3,4}$ and $c_{4,5}$. Now, the coboundary (\ref{nedj}) can be
added in order to eliminate the constant $c_{3,4}.$ The constants
are as follows:
\begin{equation}\nonumber
\small
\begin{array}{lcllcl}
c_{3,11}&\!\!=\!\!&
\frac{1}{88}((11-k)(k+2\lambda-12)\,c_{3,10}-4\,c_{4,10}),
&c_{3,4}&\!\!=\!\!&0,\\[2mm]
c_{3,13}&\!\!=\!\!&
\frac{1}{130}((13-k)(k+2\lambda-14)\,c_{3,12}-4\,c_{4,12}),&
c_{5,9}&\!\!=\!\!&
\frac{1}{10}(k-11)(k+2\lambda-12)\,c_{9,4}\\[1mm]
&&&&&+7\,c_{10,4},\\[2mm]
c_{5,10}&\!\!=\!\!&
\frac{1}{10}(k-12)(k+2\lambda-13)\,c_{10,4}+\frac{44}{5}\,c_{11,4},&
c_{5,11}&\!\!=\!\!&
\frac{1}{10}(k-13)(k+2\lambda-14)\,c_{11,4}\\[1mm]
&&&&&+\frac{54}{5}\,c_{12,4},\\[2mm]
c_{6,8}&\!\!=\!\!&
\frac{1}{18}(63\,c_{9,5}-(k-12)(k-11+2\,\lambda)\,c_{5,8} ),&
c_{6,9}&\!\!=\!\!&
\frac{1}{18}((-k+13)(k-12+2\,\lambda)\,c_{5,9})\\[1mm]
&&&&&+\frac{40}{9}\,c_{10,5},\\[2mm] c_{6,10}&\!\!=\!\!&
\frac{1}{18}(99\,c_{11,5}-(k-14)(k-13+2\,\lambda)\,c_{5,10} ),&
c_{7,8}&\!\!=\!\!&
\frac{1}{14}(35\, c_{10,6}-(k-13)\, c_{6,9}),\\[2mm]
c_{7,9}&=& \frac{1}{28}(54\, c_{9,6}-(k-12)(k-11+2\,\lambda)\,
c_{6,8}),&c_{11,4}&\!\!=\!\!&0,\\[2mm]
c_{3,12}&\!\!=\!\!&
\frac{1}{108}((12-k)(k+2\lambda-13)\,c_{3,11}-4\,c_{4,11}),&
c_{4,5}&\!\!=\!\!&1.
\end{array}
\end{equation}
Here we omit the expressions of $c_{10,4}$ and $c_{12,4}$ as they
are too long.

3) If $\lambda=\frac{1}{2}(1-k),$ then the cohomology is
two-dimensional.
\subsubsection{The case where $\mu-\lambda \geq 15$}
Let us prove that the system (\ref{ref}) has solutions that depend
on one parameter for all $\lambda.$ We have seen in the previous
section that the solutions to the system (\ref{ref}) depend on one
parameter if $\lambda$ is generic and on two parameters if
$\lambda=\frac{1}{2}(1-k)$ or $ \frac{1}{2}(1-k\pm
\sqrt{12k-23}).$ But here $k\geq 15;$ we have to study (\ref{ref})
for $\alpha=5.$ For $\alpha=5, \beta=6$ and $\gamma=7,$ the system
(\ref{ref}) has one more equation
\begin{equation}\small
\label{nak}
\begin{array}{l}
\left ( \binom{10}{5}-\binom{10}{4} \right )\, c_{10,7}- \left
(\binom{11}{5}-\binom{11}{4}\right )\, c_{11,6} + \left
(\binom{12}{6}-\binom{12}{5}\right )\, c_{12,5}+ \left
(\binom{k-11}{5}+\lambda \binom{k-11}{4}\right )\, c_{6,7}\\[2mm] -\left
(\binom{k-10}{6}+\lambda \binom{k-10}{5}\right )\, c_{5,7}+ \left
(\binom{k-9}{7}+\lambda \binom{k-9}{6}\right )\, c_{5,6} =0.
\end{array}
\end{equation}

1) For $\lambda=\frac{1}{2}(1-k),$ Eq. (\ref{nak}) implies that
the constant $c_{4,5}$ is expressed in terms of $c_{3,4}.$ Once
more we omits its explicit expression. If $k=15,$ then $c_{3,4}$
generates the system and consequently the cohomology is
one-dimensional since $\beta_{3,4}=\beta_{4,5}=0.$ If $k>15,$ then
the system (\ref{ref}) adds another condition that implies
$c_{3,4}=0.$ Therefore the cohomology is zero.

2) For $\lambda =\frac{1}{2}(1-k\pm\sqrt{12k-23}),$ we proceed as
before. The cohomology is zero.

\cqfd
\begin{remark} The study of $\mathfrak{sl}(2)$-invariant differential
operators over polynomial vector fields on $\mathbb{R},$
$\Vect_{\mathrm{P}}(\mathbb{R}),$ or over smooth vector fields on
the circle, $\Vect(\mathbb{S}^1),$ (in the case of $\mathbb{S}^1$
we express such operators in an affine coordinate) is identical
with the study of $\mathfrak{sl}(2)$-invariant differential
operators over $\Vect(\mathbb{R})$. Therefore, Theorem \ref{main}
remains true whether for $\Vect(\mathbb{S}^1)$ or
$\Vect_{\mathrm{P}}(\mathbb{R})$ since its proof is based on the
classification of $\mathfrak{ sl}(2)$-invariant differential
operators.
\end{remark}
\section{Explicit 2-cocycles for $\Vect(\mathbb{R})$ and
$\mathfrak{sl}(2)$}
The following cohomology was computed by Lecomte \cite{l}:
\begin{equation} \nonumber \mathrm
H^2(\mathfrak{sl}(2);  \cD_{\lambda,\mu})= \left\{
\begin{array}{ll}
 \bbR& \mbox{if }\,(\lambda,
\mu)=(\frac{1-k}{2},\frac{1+k}{2}),
\mbox{ and } k\in \mathbb{N}\backslash \{0\},\\
0& \mbox{otherwise.}
\end{array}
\right.
\end{equation}
The 2-cocycle that spans this cohomology is given by (here
$\omega$ is the Gelfand-Fuchs cocycle (\ref{gelfuk})):
\begin{equation}
\nonumber \Omega(X,Y,\phi\,dx^\lambda)=\omega(X,Y)\,
\phi^{(k-1)}\,dx^{\frac{1+k}{2}}.
\end{equation}
The following cohomology can be deduced from the work of
Feigin-Fuchs \cite{ff} (where
$\Vect_{\mathrm{P}}(\mathbb{\mathbb{R}})$ is the Lie algebra of
polynomial vector fields) :
\begin{equation} \label{tay}{\mathrm
H}^2 (\Vect_{\mathrm{P}}(\mathbb{\mathbb{R}}); \cD_{\lambda,\mu})=
\left\{
\begin{array}{ll}
\bbR& \mbox{if } \left \{
\begin{array}{l}
(\mu,\lambda)=(1,0),\\[2mm]
\mu-\lambda=2,3,4 \mbox{ for all }
\lambda, \\[2mm]
\mu-\lambda=7,8,9,10,11 \mbox{ for all }
\lambda, \\[2mm]
\mu-\lambda=k=12,13,14 \mbox{ but } \lambda
\mbox{ is either } \frac{1-k}{2},\\[2mm]
\mbox{ or } \frac{1-k}{2}\pm \frac{\sqrt{12k-23}}{2},\\[2mm]
\end{array}
\right.
\\
\mathbb{R}^2 & \mbox{if } \left \{
\begin{array}{l}
(\lambda,\mu)=(0,5) \mbox{ or } (-4,1),\\[2mm]
(\lambda,\mu)=\left (-\frac{5}{2}\pm \frac{\sqrt{19}}{2},
\frac{7}{2}\pm \frac{\sqrt{19}}{2}\right ),
\end{array}
\right.
\\
0&\mbox{otherwise.}
\end{array}
\right.
\end{equation}
The 2-cocycles spanning (\ref{tay}) for $k=1,2,3,4,5$ and $6$ are
as follows (here $X=f \frac{d}{dx}, Y=~g\frac{d}{dx}$):

(i) For $(\lambda,\mu)=(0,1),$ the 2-cocycle is given by
\begin{equation}
\label{ch} \Omega_{1}(X,Y,\phi\,dx^\lambda)=\omega(X,Y)\,
\phi\,dx^\lambda.
\end{equation}

(ii) For $\mu-\lambda=2,$ the 2-cocycle is given by
\begin{equation}\label{corr2}
\Omega_{2}(X,Y)=c_{1,2}\,\omega(X,Y)\, \frac{d}{dx}+c_{1,3}\,
\left |
\begin{array}{ll}
f' & g'\\
f''' & g'''
\end{array}
\right |,
\end{equation}
where  $c_{1,2}=1$ and $c_{1,3}=0$ for $\l=-\frac12;$ whereas
$c_{1,2}=0$ and $c_{1,3}=1$ for $\l\neq-\frac{1}{2}.$

(iii) For $\mu-\lambda=3,$ the 2-cocycle is given by
\begin{equation}\small \label{corr3}
\Omega_{3}(X,Y)=c_{1,2}\,\omega(X,Y)\,\frac{d^2}{dx^2}+c_{1,3}\,
\left |
\begin{array}{ll}
f' & g'\\
f^{'''} & g^{'''}
\end{array}\right |\, \frac{d}{dx}
+\frac{\l}{2}(c_{1,2}-c_{1,3})\,\left |
\begin{array}{ll}
f' & g'\\
f^{(4)} & g^{(4)}
\end{array}
\right |,
\end{equation}
where $c_{1,2}=1$ and $c_{1,3}=0$ for $\l=-1;$ whereas
$c_{1,2}=0$ and $c_{1,3}=1$ for $\l\not=-1.$

(iv) For $\mu-\lambda=4,$ the 2-cocycle is given by
\begin{equation}\small \label{corr4}
\begin{array}{ll}
\Omega_{4}(X,Y)&=c_{1,2}\,\omega(X,Y)\,
\frac{d^3}{dx^3}+\frac{1}{2}((1+2\l)c_{1,3}-(1+3\l)c_{1,2})\,
\left |
\begin{array}{ll}
f' & g'\\
f^{(4)} & g^{(4)}
\end{array}
\right |\,\frac{d}{dx}\\
&+ c_{1,3}\, \left |
\begin{array}{ll}
f' & g'\\
f^{'''} & g^{'''}
\end{array}
\right
|\,\frac{d^2}{dx^2}+\frac{\l}{10}((1-3\l)c_{1,2}+(1+2\l)c_{1,3})\,
\left |
\begin{array}{ll}
f' & g'\\
f^{(5)} & g^{(5)}
\end{array}
\right |,
\end{array}
\end{equation}
where $c_{1,3}=0$ and $c_{1,2}=1$ for $\l=-\frac{3}{2};$ whereas
$c_{1,3}=1$ and $c_{1,2}=0$ for $\l\not=-\frac{3}{2}.$

(v) For $\mu-\lambda=5,$ the two 2-cocycles are given by (where
$\alpha$ and $\beta$ are constants):
\begin{equation}\label{corr5}\small
\begin{array}{lcl}
\Omega_{5}(X,Y)&=&3 \alpha (1+\lambda)(1+2\lambda) \,\omega(X,Y)\,
\frac{d^4}{dx^4}+ 2\alpha (1+3\lambda+6\lambda^2)\left |
\begin{array}{ll}
f' & g'\\
f^{(3)} & g^{(3)}
\end{array}
\right |\,\frac{d^3}{dx^3}\\
& &+3\alpha (1+\lambda)(1+4\lambda)\left |
\begin{array}{ll}
f' & g'\\
f^{(4)} & g^{(4)}
\end{array}
\right |\,\frac{d^2}{dx^2}-\frac{1}{5}\alpha
\lambda(1+9\lambda)\left |
\begin{array}{ll}
f' & g'\\
f^{(6)} & g^{(6)}
\end{array}
\right |\\
&&+\beta\, \left |
\begin{array}{ll}
f''' & g'''\\
f^{(4)} & g^{(4)}
\end{array}
\right |.
\end{array}
\end{equation}
(vi) For $\mu-\lambda=6,$ the two 2-cocycles are given by (where
$\alpha$ and $\beta$ are constants):
\begin{equation}\label{corr6}\small
\begin{array}{ll}
\Omega_{6}(X,Y)&=\alpha (4+3\lambda(5+2\lambda))\,\omega(X,Y)\,
\frac{d^5}{dx^5}+5\alpha(2+\lambda(4+3\lambda)) \left |
\begin{array}{ll}
f'& g'\\
f^{(3)}& g^{(3)} \end{array} \right |\,
\frac{d^4}{dx^4}\\
&+5\alpha(\lambda(3+4\lambda)-2) \left |
\begin{array}{ll}
f'& g'\\
f^{(4)}& g^{(4)} \end{array} \right |\,
\frac{d^3}{dx^3}+5\alpha(2+\lambda(4+3\lambda)) \left |
\begin{array}{ll}
f'& g'\\
f^{(5)}& g^{(5)} \end{array} \right |\,\frac{d^2}{dx^2}\\
&+\beta
\left |\!
\begin{array}{ll}
f^{(3)}& g^{(3)}\\
f^{(4)}& g^{(4)} \end{array} \!\!\!\right |\,
\frac{d}{dx}+\alpha(4+15\lambda+6\lambda^2) \left |\!
\begin{array}{ll}
f'& g'\\
f^{(6)}& g^{(6)} \end{array} \!\!\!\right |\,
\frac{d}{dx}-\frac{\lambda}{5}\beta\, \left |\!
\begin{array}{ll}
f^{(3)}& g^{(3)}\\
f^{(5)}& g^{(5)} \end{array} \!\!\!\right |.
\end{array}
\end{equation}
In order to complete the list of 2-cocycles spanning (\ref{tay}) we
need the following two Lemmas.
\begin{lemma}
\label{baid} Every $2$-cocycle in
$\mathrm{H}^2(\Vect_{\mathrm{P}}(\mathbb{R}); {\cal
D}_{\lambda,\mu}) $ can be reduced to a $2$-cocycle vanishing on
$\mathfrak{sl}(2),$ except those given in
$(\ref{ch})-(\ref{corr6})$.
\end{lemma}
\begin{proof} Consider a general form of a 2-cocycle
(where $X=f\frac{d}{dx}, Y=g\frac{d}{dx}\in
\Vect_{\mathrm{P}}(\mathbb{R})$ and $\phi\,dx^\lambda\in
\cF_\lambda$):
\begin{equation}
\label{hab} c(X,Y,\phi\,dx^{\lambda})=\sum_{i+j+l=k+2}\,c_{i,j}\,
f^{(i)}\,g^{(j)}\,\phi^{(l)}dx^{\lambda+k}.
\end{equation}
We will eliminate coboundaries in order to turn the 2-cocycle
above into a 2-cocycle vanishing on $\mathfrak{sl}(2).$ Consider a
general expression of a coboundary
\[
\begin{array}{ccl}
\delta B(X,Y,\phi\,dx^\lambda)&=&-\beta_0\,fg'\,
\phi^{(k+1)}-\beta_0(\binom{k+1}{\alpha}+\lambda
\binom{k+1}{\alpha-1})\, fg^{(\alpha)}\,
\phi^{(k+2-\alpha)}\\[2mm]
&&-\sum_{\alpha\geq2} \beta_1(\binom{k}{\alpha}+\lambda
\binom{k}{\alpha-1}) f'g^{(\alpha)}\, \phi^{(k+1-\alpha)}+
\mbox{higher order terms}\\[2mm]
&&-(f\leftrightarrow g).
\end{array}
\]
Immediately we see that the constant $c_{0,1}$ can be eliminated
upon putting $c_{0,1}=-\beta_0.$ On the other hand, the 2-cocycle
condition implies that $c_{\gamma,0}=-c_{0,1}\,
(\binom{k+1}{\gamma}+\lambda\, \binom{k+1}{\gamma-1} ) .$

1) For $k=1,$ the 2-cocycle takes the form
\begin{equation}
\nonumber \Omega_{1}(X,Y,\phi)=c_{1,2}\,\omega(X,Y)\, \phi.
\end{equation}
On the other hand, the coboundary takes the form
\[
\delta B(X,Y,\phi)=\lambda\,\alpha_1\,\omega(X,Y)\, \phi,
\]
where $\alpha_1$ is a constant. The 2-cocycle is trivial except
for $\lambda=0.$\\

\noindent 2) For $k=2,3,4,5,6,$ we proceed as before.

Suppose now that $k>6.$ We will deal with the coefficients
$c_{1,\gamma}.$ The 2-cocycle condition implies that the component
of $f'\, g^{\beta}\, h^{\gamma}\, \phi^{k+2-\beta-\gamma}$, which
should be zero, is equal to
\begin{equation}
\label{kh}
\begin{array}{l}
c_{\beta+\gamma-1,1}\left (\binom{\beta+\gamma-1}
{\beta}-\binom{\beta+\gamma-1}{\beta-1}\right ) -c_{1,\gamma}\left
(\binom{k+1-\gamma}{\beta}+\lambda
\binom{k+1-\gamma}{\beta-1}\right )\\[2mm]
+c_{1,\beta}\left (\binom{k+1-\beta}{\gamma}+\lambda
\binom{k+1-\beta}{\gamma-1}\right )=0.
\end{array}
\end{equation}
We have two cases:

i) For $\lambda=\frac{1-k}{2}.$ In this case, the coefficient of
$f'g''\, \phi^{k-1}$ is zero in the expression of the coboundary.
But $c_{1,3}$ can be eliminated upon putting $
c_{1,3}=\frac{1}{6}\, k(k-1) \left (k-2+3\lambda\right
)\,\beta_1.$ By putting $\beta=2,$ we can see from 
(\ref{kh}) that all $c_{t,1}$ can be expressed in terms of
$c_{1,2}$. They are given by the induction formula:
\begin{equation}
\label{upon}
\begin{array}{l}
c_{1,i}=\frac{2}{i-3}\, \left (-c_{1,i-1}\left
(\binom{k+2-i}{2}+\lambda \,\binom{k+2-i}{1}\right
)+c_{1,2}(\binom{k-1}{i-1}+\lambda \,\binom{k-1}{i-2})\right
)\quad \mbox{for } i>3.
\end{array}
\end{equation}
However, for $\beta=3$ and $\gamma=4$ the system (\ref{kh})
becomes
\[
\binom{k-1}{4}\, (1+k)(1+3k)\,c_{1,2}=0.
\]
As $k>4,$ the equation above admits a solution only for
$c_{1,2}=0.$ Thus, all $c_{1,\gamma}$ are zero.

ii) If $\lambda\not =\frac{1-k}{2},$ then the constant $c_{1,2}$
can be eliminated and we proceed as before.

Now we deal with the coefficients $c_{2,s}.$ These coefficients
can be eliminated upon taking
\[
\begin{array}{ccl}
\beta_{s+1,k-s-1}&=&\frac{1}{(s+1)(s-2)}\left( c_{2,s} +
(k-s)(k-s-1+2\,\lambda)\, \beta_{s,k-s}\right )\\[2mm]
&&+\frac{1}{(s+1)(s-2)}\left(-2\, \left (\binom{k-2}{s}+\lambda
\binom{k-2}{s-1}\right )\, \beta_{2,k-2}\right ). \end{array} \]
Finally, the remaining 2-cocycle vanishes on $\mathfrak{sl}(2).$
\end{proof}
\begin{lemma}
\label{saa} Every coboundary $\delta  (B)\in
B^2(\Vect(\mathbb{R}); \cD_{\lambda,\mu})$ vanishing on
$\mathfrak{sl}(2)$ possesses the following properties. The
operator $B$ coincides (up to a nonzero factor) with the
transvectant $J_{k+1}^{-1,\lambda}$, where
$\gamma_{0,k+1}=\gamma_{1,k}=0.$ In addition (here
$X=f\frac{d}{dx}, Y=g\frac{d}{dx}\in \Vect(\mathbb{R})$ and
$\phi\, dx^{\lambda}\in \cF_\lambda$)
\begin{equation}
\label{nedj2}\delta (B)(X,Y,\phi\,dx^\lambda)=
\sum_{i+j+l=k+2}\,\beta_{i,j}\,f^{(i)}\,g^{(j)}\,\phi^{(l)}\,
dx^{\lambda+k},
\end{equation}
where
$$
\beta_{0,j}=\beta_{1,j}=\beta_{2,j}=0,
$$
and \small{
$$
\begin{array}{lcl}
\beta_{3,4}&=&\frac{1}{24}\binom{k-2}{3}\left(k^2 + 4 (\lambda
-1)\lambda + k (4\lambda-5)\right
)((k-1)(k-2+3\lambda)\gamma_{2,k-1}-( k -1+ 2 \lambda
)\gamma_{3,k-2})
\\[3mm]
\beta_{4,5}&=&-\frac{1}{480}\binom{k-2}{5} (k-3+ 2\lambda )(k^3 +
4(\lambda-1) \lambda (2\lambda-19) + 3 k^2(2
\lambda-7 )+2k(49+6(\lambda-7 )\lambda))\\[3mm]
&&\times ( (k-1)(k-2+3\lambda)\gamma_{2,k-1} - (k -1+ 2 \lambda
)\gamma_{3,k-3}).
\end{array}
$$
}
\end{lemma}
\begin{proof}
Similar to Proposition \ref{maa}.
\end{proof}

Now we will explain how we can deduce the explicit expressions of
the 2-cocycles that span
$\mathrm{H}^2(\Vect_{\mathrm{P}}(\mathbb{R});
\cal{D}_{\lambda,\mu})$ by using the results of Sec. \ref{hop}. To
save space, we give details of the computation only for
$\mu-\lambda=7,8,9,10,11.$ The other cases, namely
$\mu-\lambda=12,13,14,$ can be deduced by the same way. We start
with any 2-cocycle $c\in Z^2(\Vect_{\mathrm{P}}(\mathbb{R});
\cal{D}_{\lambda,\mu})$ vanishing on $\mathfrak{sl}(2).$ This is
actually possible, thanks to Lemma \ref{baid}. The 2-cocycle
condition of $c$ has already been studied Sec. \ref{coc}. The
2-cocycle $c$ is generated by the two constants $c_{3,4}$ and
$c_{4,5}.$ We have the following cases:

1) $\lambda=\frac{1-k}{2}.$ By Lemma \ref{saa}, one of the
constants $c_{3,4}$ or $c_{4,5}$ can be eliminated by adding a
coboundary with an appropriate value of $\gamma_{2,k-2}.$ We
obtain, therefore, a unique 2-cocycle that is non-trivial in
$\mathrm{H}^2(\Vect_{\mathrm{P}}(\mathbb{R});
\cal{D}_{\lambda,\mu}).$

2) $\lambda =\frac{2-k}{3}.$ By Lemma \ref{saa}, one of the
constants $c_{3,4}$ or $c_{4,5}$ can be eliminated by adding a
coboundary with an appropriate value of $\gamma_{3,k-3}.$ We
obtain, therefore, a unique 2-cocycle that is non-trivial in
$\mathrm{H}^2(\Vect_{\mathrm{P}}(\mathbb{R});
\cal{D}_{\lambda,\mu}).$

3) $\lambda$ is a solution to the equation
$$
k^2 + 4 (\lambda -1)\lambda + k (4\lambda-5)=0. $$Then
$\beta_{3,4}=0.$ Therefore the constant $c_{4,5}$ can be
eliminated with an appropriate value of $\gamma_{2,k-1}.$ We
obtain, therefore, a unique 2-cocycle that is non-trivial in
$\mathrm{H}^2(\Vect_{\mathrm{P}}(\mathbb{R});
\cal{D}_{\lambda,\mu}).$\\

4) $\lambda$ is a solution to the equation
$$(k-3+ 2\lambda )(k^3 +
4(\lambda-1) \lambda (2\lambda-19) + 3 k^2(2 \lambda-7
)+2k(49+6(\lambda-7 )\lambda))=0.
$$
Then $\beta_{4,5}=0.$ Therefore the constant $c_{3,4}$ can be
eliminated with an appropriate value of $\gamma_{2,k-1}.$ We
obtain, therefore, a unique 2-cocycle that is non-trivial in
$\mathrm{H}^2(\Vect_{\mathrm{P}}(\mathbb{R});
\cal{D}_{\lambda,\mu}).$

5) If $\lambda$ is not as in 1)--4). Whatever the value of
$\lambda$ is the constant $c_{3,4}$ can be eliminated with an
appropriate value of $\gamma_{2,k-1}.$ We obtain, therefore, a
unique 2-cocycle that is non-trivial in
$\mathrm{H}^2(\Vect_{\mathrm{P}}(\mathbb{R});
\cal{D}_{\lambda,\mu}).$
\subsection{Further remarks}
It would be interesting to study the cohomology arising in the
deformation of symbols at the group level, $\Diff(\mathbb{R}).$ We
do not know whether our 2-cocycles introduced here can be
integrated to the group. Nevertheless,  the 2-cocycle (\ref{ch})
can be integrated to a 2-cocycle $A\in
\mathrm{H}^2(\Diff(\mathbb{R}); \cD_{\lambda,\lambda+1})$ (here
$F, G\in \Diff(\mathbb{R})$ and $\phi\, dx^{\lambda}\in
\cF_\lambda$):
$$
A(F,G,\phi\,dx^\lambda):=\log(F\circ G)'\, \frac{G''}{G'} \,
\phi\,dx^{\lambda+1}.
$$
This 2-cocycle is just the multiplication operator by the
well-know  Bott-Thurston cocycle \cite{bott}. Let
$S(f):=\frac{f'''}{f'}-\frac{3}{2} \left (\frac{f''}{f'}\right
)^2$ be the Schwarz derivative. Then the 2-cocycle (\ref{or}) can
be integrated to $B\in \mathrm{H}^2(\Diff(\mathbb{R}),
\PSL(2,\mathbb{R}); \cD_{\lambda,\lambda+5}):$
$$
B(F,G,\phi\,dx^\lambda):=\left |
\begin{array}{ll}
G^* S(F) & S(F)\\
G^* S(F)' & S(F)'
\end{array}
\right | \, \phi\,dx^{\lambda+5}.
$$
This 2-cocycle is also the multiplication operator by a 2-cocycle
introduced by Ovsienko-Roger \cite{or2}.

It would also be interesting to study the cohomology arising in
the context of deformation of the space of symbols on
multi-dimensional manifolds.\\

\noindent {\bf Acknowledgments.} I would like to thank 
M. Ben Ammar, D. Leites, V. Ovsienko and J. Stasheff for their 
suggestions and remarks.


\begin{thebibliography}{99}
\bibitem{aalo} Agrebaoui B, Ammar F, Lecomte P and Ovsienko V,
{\rm Multi-parameter deformations of the module of symbols of
differential operators.} {\it Int. Math. Res. Not.,} no. 16,
(2002), 847--869.

\bibitem{bott} Bott R, {\rm
On the charachteristic classes of groups of diffeomorphisms}, {\it
Enseign. Math.} {\bf 23:3-4,} (1977), 209--220.

\bibitem{b1} Bouarroudj S, {\rm
Cohomology of the vector fields Lie algebras on $\mathbb{RP}^1$
acting on bilinear differential operators}, {\it Int. Jour. Geom.
Methods. Mod. Phys.} {\bf 2,} no. 1, (2005), 23--40.

\bibitem{b2} Bouarroudj S, {\rm Projective and conformal Schwarzian
derivatives and cohomology of Lie algebras vector fields related
to differential operators,} {\it Int. Jour. Geom. Methods. Mod.
Phys.} {\bf 3,} no. 4, (2006), 667-696.

\bibitem{bo} Bouarroudj S and Ovsienko V, {\rm Three cocycles on
$\Diff(S^1)$ generalizing the Schwarzian derivative}, {\it
Internat. Math. Res. Notices,} no. 1,1998, 25--39.

\bibitem{do}
Duval C and Ovsienko V, {\rm Conformally equivariant quantum
Hamiltonians.} {\it Selecta Math. (N.S.)} {\bf 7}, no. 3, (2001),
291--320.

\bibitem{ff}
Feigin B L and Fuchs D B,
{\rm Homology of the Lie algebra of vector fields on the line},
{\it Func. Anal. Appl.,} 14 (1980), 201--212.

\bibitem{f} Fuchs D B, {\it Cohomology of infinite-dimensional
Lie algebras}, Contemp. Soviet. Math., Consultants Bureau,
New-York, 1986.

\bibitem{gar}
Gargoubi H, {\rm Sur la g\'eom\'etrie de l'espace des
op\'erateurs diff\'erentiels lin\'eaires sur $\mathbb{R}.$} {\it
Bull. Soc. Roy. Sci. Li\`ege.} Vol. 69, 1, (2000), 21--47.

\bibitem{gf}
Gelfand G F and Fuchs D B, {\rm Cohomology of the Lie algebra
of vector fields on the circle}, {\it Func. Anal. Appl.,} {\bf
2:4}, (1968), 342--343.


\bibitem{g} Gordan P,
{\it Invariantentheorie,} Teubner, Leipzig, 1887.

\bibitem{l}
Lecomte P B A, {\rm On the cohomology of
$\mathfrak{sl}(n+1,\mathbb{R})$ acting on differential operators
and $\mathfrak{sl}(n+1,\mathbb{R})$-equivariant symbols,} {\it
Indag. Math. NS.} 11 (1), (2000), 95--114.

\bibitem{lo}
Lecomte P B A and Ovsienko V, {\rm Cohomology of the vector
fields Lie algebra and modules of differential operators on a
smooth manifold,} {\it Compositio Mathematica.} {\bf 124:} no.1,
2000, 95--110.

\bibitem{los}
Losik M V, {\rm Cohomology of the Lie algebra of vector fields
with nontrivial coefficients}, {\it Func. Anal. Appl.,} {\bf 6,}
(1972), 289--291.

\bibitem{nr}
Nijenhuis A and Richardson R W, Deformation of homomorphisms of
Lie algebras, Bull. AMS, {\bf 73} (1967) 175--179.

\bibitem{olv}
Olver P, Applications of Lie groups to differential
equations. {\rm Springer,} 1993.

\bibitem{or2} Ovsienko V and Roger C, {\rm
Generalization of Virasoro group and Virasoro algebra through
extensions by modules of tensor-densities on $S^1$.} {\it Indag.
Math., (N.S.),} {\bf 9}, no.2, (1998), 277--288.

\bibitem{ot}
Ovsienko V and Tabachnikov S, {\it Projective differential
geometry old and new: from the Schwarzian derivative to cohomology
of diffeomorphism groups,} {\rm Cambridge University Press,} 2004.

\bibitem{tsu} Tsujishita T, {\rm On the continuous cohomology of
the Lie algebra of vector fields,} {\it Proc. Japan Acad.}, A 53,
N. 4, (1977), 134--138.

\end{thebibliography}
\end{document}